\documentclass[10pt,journal]{IEEEtran}

\usepackage{xspace,amssymb,epsfig,subfigure,syntonly}
\usepackage{epsfig,amsmath,color}
\usepackage{xspace,syntonly}
\usepackage{url}
\usepackage{algorithm,algorithmic}

\newcommand{\utwi}[1]{\mbox{\boldmath $#1$}}
\newcommand{\E}{{\mathbb{E}}}

\newcommand{\Prob}{{\textrm{Pr}}}
\newcommand{\trace}{{\textrm{Tr}}}
\newcommand{\rank}{{\textrm{rank}}}

\newcommand{\cD}{{\cal D}}

\newcommand{\cN}{{\cal N}}

\newcommand{\cA}{{\cal A}}

\newcommand{\cT}{{\cal T}}

\newcommand{\cE}{{\cal E}}

\newcommand{\cF}{{\cal F}}

\newcommand{\cU}{{\cal U}}

\newcommand{\cH}{{\cal H}}

\newcommand{\bd}{{\bf d}}
\newcommand{\be}{{\bf e}}

\newcommand{\bp}{{\bf p}}
\newcommand{\bq}{{\bf q}}

\newcommand{\bu}{{\bf u}}
\newcommand{\bv}{{\bf v}}
\newcommand{\bi}{{\bf i}}
\newcommand{\bz}{{\bf z}}
\newcommand{\by}{{\bf y}}
\newcommand{\bA}{{\bf A}}
\newcommand{\bB}{{\bf B}}

\newcommand{\bM}{{\bf M}}

\newcommand{\bS}{{\bf S}}

\newcommand{\bI}{{\bf I}}

\newcommand{\bY}{{\bf Y}}
\newcommand{\bV}{{\bf V}}

\newcommand{\bPi}{{\utwi{\Pi}}}

\newcommand{\sfH}{\textsf{H}}
\newcommand{\sfT}{\textsf{T}}

\usepackage{epstopdf}
\usepackage[noadjust]{cite}

\begin{document}

\newtheorem{definition}{Definition}
\newtheorem{remark}{Remark}
\newtheorem{proposition}{Proposition}
\newtheorem{lemma}{Lemma}

\title{Optimal Dispatch of Residential Photovoltaic Inverters Under Forecasting Uncertainties} 
\author{Emiliano Dall'Anese, \emph{Member}, \emph{IEEE}, Sairaj V. Dhople, \emph{Member}, \emph{IEEE}, Brian B. Johnson, \emph{Member}, \emph{IEEE},  \\ and Georgios B. Giannakis, \emph{Fellow}, \emph{IEEE}
\thanks{\protect\rule{0pt}{0.5cm}%

Paper submitted on June 29, 2014; revised September 21, 2014; accepted October 16, 2014. This work was supported by the Institute of Renewable Energy and the Environment (IREE) grant no. RL-0010-13, University of Minnesota, by the Laboratory Directed Research and Development (LDRD) Program at the National Renewable Energy Laboratory, and by the National Science Foundation (NSF) grants CCF 1423316 and CyberSEES 1442686. 

E. Dall'Anese, S. V. Dhople, and G. B. Giannakis are with the  Department of Electrical and Computer Engineering, and also the Digital Technology Center, University of Minnesota, 200 Union Street SE, Minneapolis, MN, USA; e-mail: {\tt \{emiliano, sdhople, georgios\}@umn.edu}. B. B. Johnson is with the National Renewable Energy Laboratory, Golden, CO, USA; e-mail: {\tt brian.johnson@nrel.gov}.   
}
}

\markboth{IEEE JOURNAL OF PHOTOVOLTAICS}
{Dall'Anese \MakeLowercase{\textit{et al.}}: }

\maketitle

\vspace{.5cm}

\begin{abstract}

Efforts to ensure reliable operation of existing low-voltage distribution systems with high photovoltaic (PV) generation have focused on the possibility of inverters providing ancillary services such as active power curtailment and reactive power compensation. Major benefits include the possibility of averting overvoltages, which may otherwise be experienced when PV generation exceeds the demand. This paper deals with ancillary service procurement in the face of solar irradiance forecasting errors. In particular, assuming that the forecasted PV irradiance can be described by a random variable with known (empirical) distribution, the proposed uncertainty-aware optimal inverter dispatch (OID) framework indicates which inverters should provide ancillary services with a guaranteed a-priori risk level of PV generation surplus. To capture forecasting errors, and strike a balance between risk of overvoltages and (re)active power reserves, the concept of conditional value-at-risk is advocated. 
Due to AC power balance equations and binary inverter selection variables, the formulated OID involves the solution of a nonconvex mixed-integer nonlinear program. However, a computationally-affordable convex relaxation is derived by leveraging  sparsity-promoting regularization approaches and semidefinite relaxation techniques. 
\end{abstract}


\section{Introduction}
\label{sec:Introduction}

Deployment of photovoltaic (PV) systems in residential settings
promises a multitude of environmental and economic advantages,
including a sustainable capacity expansion of distribution
systems. However, a unique set of challenges related to power quality,
efficiency, and reliability may emerge, especially when an increased
number of PV systems are deployed in existing distribution networks,
and operate according to current practices~\cite{Liu08,Cpuc12}. One
challenge is associated with overvoltages when PV generation exceeds demand~\cite{Tonkoski11}. 

To ensure reliable operation of existing distribution feeders even
during peak PV generation hours, recent efforts have focused on the
possibility of inverters providing ancillary services~\cite{Tricoli10,
  Caramanis_PESTD14,Craciun13}.  For instance, reactive power
compensation approaches have been recognized as a viable option to
effect voltage regulation at the medium-voltage distribution
level~\cite{Turitsyn11,Aliprantis13,Bolognani13}. The amount of
reactive power injected or absorbed by inverters can be computed based
on either local droop-type proportional
laws~\cite{Turitsyn11,Aliprantis13} or optimal power flow (OPF)
strategies~\cite{Bolognani13}. Either way, voltage regulation with
this approach comes at the expense of low power factors at the
substation and high network currents, with the latter leading to high
power losses in the network~\cite{Tonkoski11}. Alternative approaches
rely on operating inverters at unity power factor while curtailing
part of the available active power~\cite{Tonkoski11,Samadi14}. Active
power curtailment strategies are particularly effective in the
low-voltage portion of distribution feeders, where the high
resistance-to-inductance ratio of low-voltage overhead lines renders
voltage magnitudes more sensitive to variations in the active power
injections.
An optimal inverter dispatch (OID) framework was proposed in~\cite{OID} to set both active and reactive power setpoints so that the network operation is optimized according to chosen criteria (e.g., minimizing power losses), while ensuring voltage regulation and adhering to other electrical constraints.

The approaches
in~\cite{Turitsyn11,Aliprantis13,Bolognani13,Tonkoski11,Samadi14,OID,OID_TEC}
are suitable for \emph{real-time} network control, where the setpoints
of the inverters scheduled to provide ancillary services are
fine-tuned based on instantaneous load measurements and prevailing
ambient conditions. Distinct
  from~\cite{Turitsyn11,Aliprantis13,Bolognani13,Tonkoski11,Samadi14,OID,OID_TEC},
  the problem of ancillary service \emph{procurement} is considered in
  this paper. Specifically, ancillary service procurement refers here
  to the task of scheduling the inverters that will be required to
  provide ancillary services (in e.g., minute-, hour- or day-ahead
  markets~\cite{Caramanis_PESTD14,Craciun13}), as well as quantifying
  both reactive reserves of the selected inverters and the active
  powers that inverters may be required to curtail. In this case,
  system operators cannot solely rely on the expected irradiance
  conditions to quantify the amount of ancillary services to
  provision, and irradiance forecasting errors must be taken into
  account. In fact, an excess of generation (compared to the expected
  one) may require additional inverters other than the ones scheduled (without accounting for forecasting errors) to provide ancillary
  services in order to avoid overvoltages.

The OID framework recently
proposed in~\cite{OID} is considerably broadened here by leveraging tools from risk-aware portfolio optimization to account for solar irradiance forecasting errors. Distinct from the real-time optimization method in~\cite{OID}, the approach developed in this paper enables effective provisioning of ancillary services in minute-, hour-, and day-ahead markets~\cite{Caramanis_PESTD14,Craciun13}, by  \emph{identifying} the subset of critical PV inverters that will strongly impact both voltages and network performance objectives, and \emph{quantifying} the amount of ancillary services that should be secured from each of the selected PV inverters. Specifically, for a given distribution of forecasting errors, the novel uncertainty-aware OID returns the amount of active power that can be curtailed in order to ensure voltage regulation with arbitrarily high probability and the amount of reactive power necessary to fulfill additional objectives.  The proposed scheme is grounded on an AC power flow model, and it involves the solution of an OPF type problem encapsulating well-defined performance criteria and operational constraints. To capture forecasting errors, the conditional value-at-risk (CVaR) is advocated~\cite{McNeilbook,Rockafellar00}, and utilized to trade off risks of overvoltage conditions for active power curtailment and reactive power compensation capabilities. 
Further, the resultant uncertainty-aware OID scheme involves the solution of a \emph{convex} program and handles arbitrary probability distributions for the forecasting errors~\cite{Bacher2009,Lorenz09}. 

Related prior works include, e.g.,~\cite{Sjodin12}, where chance-constrained optimal power flow (OPF) approaches were considered for high-voltage transmission systems with uncertain wind generation; a DC power flow approximation was utilized, along with Gaussian-distributed wind forecasting errors. Multi-period DC OPF was considered in e.g.,~\cite{Warrington13,Summers14}, where generation uncertainty was accounted for, while computing the schedule for controllable devices that minimize the expected operating costs. Upper bounds on the chance constraints based on e.g., Markov and Chebyshev inequalities were explored in~\cite{Summers14}. Finally, extensions to the unit commitment problem can be found in~\cite{Wang12}.  At the distribution level, a chance-constrained DC OPF formulation was developed in~\cite{Chertkov13} to mitigate the effects of  Gaussian-distributed forecasting errors on line currents and voltages. However, the DC power flow approximation may not be suitable for low-voltage resistive networks. Additional distributions for the renewable generation were considered in~\cite{Kloppel13}, where a nonconvex chance-constrained AC OPF was formulated, and solved via off-the-shelf routines for nonlinear (nonconvex) programs. An economic dispatch problem in the presence of uncertain wind generation was proposed in~\cite{YuSGComm13}; in lieu of chance constraints, the cost of the problem was regularized with CVaR-type terms capturing the risk of generation shortage. 

The remainder of the paper is organized as
follows. System modeling is outlined in Section~\ref{sec:Formulation}, along with an overview of the OID with perfect knowledge of solar irradiance~\cite{OID}. Basics of CVaR are provided in Section~\ref{sec:CVaR}, whereas the uncertainty-aware-OID is outlined in Section~\ref{sec:robustOID}. Case studies are discussed in~\ref{sec:Simulations}, while Section~\ref{sec:Conclusions} concludes the paper.\footnote{Notation:
  Upper-case (lower-case) boldface letters will be used for matrices
  (column vectors); $(\cdot)^\sfT$ for transposition; $(\cdot)^*$
  complex-conjugate; and, $(\cdot)^\sfH$ complex-conjugate
  transposition; $\Re\{\cdot\}$ and $\Im\{\cdot\}$ denote the real and
  imaginary parts of a complex number, respectively; $\mathrm{j} :=
  \sqrt{-1}$ the imaginary unit. $\mathbb{R}_+ := \{x \in \mathbb{R}: x\geq 0\}$; 
  $\trace(\cdot)$ the matrix trace;
  $\rank(\cdot)$ the matrix rank; $|\cdot|$ denotes the magnitude of a
  number or the cardinality of a set; $\|\bv\|_2 := \sqrt{\bv^\sfH
    \bv}$; and $\|\cdot\|_F$ stands
  for the Frobenius norm. For any $x \in \mathbb{R}$, $[x]_+ := \max\{0,x\}$.
  $\mathbb{I}_{\{A\}}$ is the indicator function (i.e., $\mathbb{I}_{\{A\}} = 1$ if event $A$ is true, 
  and 0 otherwise). Finally, $\bI_N$ denotes the $N \times N$
  identity matrix; and, $\mathbf{0}_{M}$, $\mathbf{1}_{M}$ the $M \times 1$ vectors with all zeroes and ones,
  respectively.}

\section{Preliminaries} \label{sec:Formulation}

\subsection{Network and inverter models}
\label{sec:models}

Consider a distribution system comprising $N+1$ nodes collected in the
set $\cN := \{0,1,\ldots,N\}$ (node $0$ denotes the secondary
of the step-down transformer), and lines represented by the set of
edges $\cE := \{(m,n)\} \subset \cN \times \cN$. Subsets
$\cU, \cH \subset \cN$ collect nodes corresponding to utility poles and households with installed
 PV inverters, respectively. For simplicity of exposition, a balanced system is
considered; however, the framework proposed subsequently can be 
extended to unbalanced multi-phase systems 
following the method in~\cite{Dallanese-TSG13}. 

Let $V_n\in \mathbb{C}$ and $I_n \in \mathbb{C}$ denote the phasors
for the line-to-ground voltage and the current injected at node $n \in
\cN$, respectively, and define $\bi := [I_0, I_1, \ldots, I_N]^\sfT \in
\mathbb{C}^{N+1}$ and $\bv := [V_0, V_1, \ldots, V_N]^\sfT \in
\mathbb{C}^{N+1}$. Using Ohm's and Kirchhoff's circuit laws, the linear
relationship $\bi = \bY \bv$ can be established, where the system
admittance matrix  $\bY \in \mathbb{C}^{N+1 \times N+1}$ is formed
based on the system topology and the $\pi$-equivalent circuits of the
lines $(m,n) \in \cE$; see e.g.,~\cite{LavaeiLow, Dallanese-TSG13,OID}.  
A constant $PQ$ model~\cite{kerstingbook} is adopted for the load,
with  $P_{\ell,h}$ and $Q_{\ell,h}$ denoting the active and reactive
demands at node $h \in \cH$, respectively (clearly, $P_{\ell,h} = Q_{\ell,h} = 0$ for all 
$h \in \cU$). 

\begin{figure}[t]
\begin{center}
\subfigure[]{\includegraphics{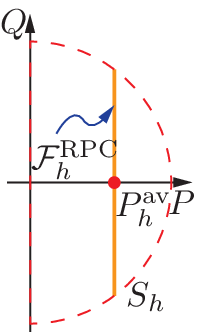}} 
\subfigure[]{\includegraphics{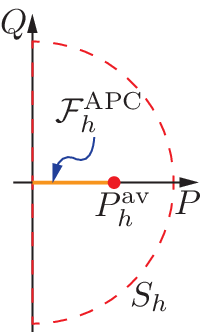}}
\subfigure[]{\includegraphics{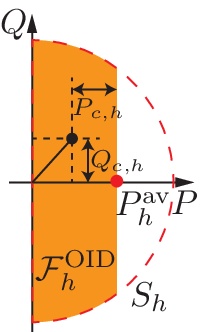}}
\subfigure[]{\includegraphics{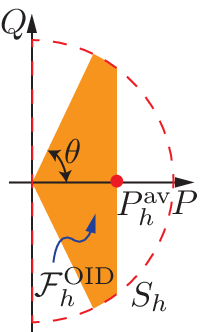}}
\end{center}
\vspace{-.4cm}
\caption{Feasible operating regions for the $h$-th inverter with apparent
  power rating $S_h$ under a) reactive power control, b) active power curtailment, c) OID with joint control of real and reactive power, and d) OID with a lower-bound on power factor~\cite{OID}.}
\label{Fig:OIDregions}
\vspace{-.4cm}
\end{figure}

For \emph{given} solar irradiation
conditions, let $P_h^{\textrm{av}}$ denote the
\emph{available active power} from the PV array at node $h \in
\cH$.  Following business-as-usual practices~\cite{Cpuc12}, grid-tied inverters  
operate at the unity-power-factor setpoint $(P_h^{\textrm{av}},0)$. 
To address emerging overvoltage and power quality concerns~\cite{Liu08},
inverters may be called upon to provide ancillary services~\cite{Caramanis_PESTD14,Tricoli10}. 
These include e.g., Volt/VAR support~\cite{Turitsyn11,Aliprantis13,Bolognani13}  
and active power curtailment~\cite{Tonkoski11}, with the allowed inverter operating 
regime on the complex-power plane 
illustrated in Fig.~\ref{Fig:OIDregions}(a) and~\ref{Fig:OIDregions}(b), respectively.  
The OID framework in~\cite{OID} offers increased flexibility over Volt/VAR support and active power
curtailment, by invoking a joint control of real and reactive powers produced by PV inverters. In
particular, the allowed operating regime 
for the PV inverter at household $h$ is illustrated in Fig.~\ref{Fig:OIDregions}(d) and described by
\begin{equation}
\label{mg-PV} 
\cF^\mathrm{OID}_h(P_{h}^{\textrm{av}}  ):= \left\{ P_{c,h}, Q_{c,h}:  \hspace{-.2cm}
\begin{array}{l}
0  \leq P_{c,h}  \leq  P_{h}^{\textrm{av}}   \\
Q_{c,h}^2  \leq  S_{h}^2 - (P_{h}^{\textrm{av}}  - P_{c,h})^2 \\
| Q_{c,h} | \leq \tan \theta (P_{h}^{\textrm{av}}  - P_{c,h}) 
\end{array}
\hspace{-.2cm} \right\} \nonumber 
\end{equation}
where $P_{c,h}$ is the active power curtailed, $Q_{c,h}$ is the
reactive power injected ($Q_{c,h} > 0$) or absorbed ($Q_{c,h} < 0$), and $S_h$ is the apparent power rating. 
In the absence of minimum power factor constraints (i.e., $\theta=\pi/2$),
the operating region corresponds to the one in Fig.~\ref{Fig:OIDregions}(c).

\subsection{Optimal inverter dispatch with known available powers}
\label{sec:Dispatch}

An overview of the OID with \emph{perfect} knowledge of $\{P_h^{\textrm{av}}\}_{h \in \cH}$  
is provided next, to lay the foundation for the uncertainty-aware framework outlined in Section~\ref{sec:Robust}.  

For given available powers $\{P_h^{\textrm{av}}\}_{h \in \cH}$, the objective of the OID is to identify the critical inverters that should be \emph{dispatched} in order to ensure electrical network constraints, 
and compute their optimal steady-state active/reactive power setpoints. To this end, let $z_h$ be a binary optimization variable indicating whether PV
inverter $h$ provides ancillary services ($z_h = 1$) or not ($z_h = 0$), and let $\bp_{c} \in \mathbb{R}_{+}^{|\cH|}$, $\bq_{c} \in \mathbb{R}^{|\cH|}$ be vectors collecting the active powers curtailed and the reactive powers injected/absorbed by the inverters. Further, let $C(\bv, \bp_{c})$ be a given cost function capturing  network- and customer-oriented objectives~\cite{OID}; for instance, $C(\bv, \bp_{c})$ may account for  active power losses in the network and possible costs associated with active power set points~\cite{Caramanis_PESTD14}.  With these definitions, a rendition of the OID task is formulated as:
\begin{subequations} 
\label{Poidnonconvex}
\begin{align} 
\hspace{1.8cm} & \hspace{-1.7cm} \min_{\bv, \bi, \bp_{c}, \bq_{c},\{z_h\}} \,\,  C(\bv, \bp_{c})  + c_z \sum_{h \in \cH} z_h \label{nonc-cost} \\
\mathrm{subject\,to} \,\, & \bi = \bY \bv, \, \{z_h\} \in \{0,1\}^{|\cH|}  \mathrm{~,and}  \nonumber  \\ 
V_h I_h^* & = (P_{h}^{\textrm{av}} - P_{c,h} - P_{\ell,h} ) + \mathrm{j}(Q_{c,h} - Q_{\ell,h}) \hspace{-.2cm} \label{nonc-balance} \\ 
 V_n I_n^* &  = 0    \hspace{3.2cm}  \forall \, n \in \cU \label{nonc-balancepoles}  \\
V^{\mathrm{min}} & \leq |V_n| \leq V^{\mathrm{max}}  \hspace{1.55cm}  \forall \, n \in \cN   \label{nonc-Vlimits} \\
(P_{c,h}, Q_{c,h}) & \in  \left\{
\begin{array}{l}
\{(0,0)\} , \mathrm{~if~} z_h = 0 \\
 \cF^\mathrm{OID}_h ,  \mathrm{~~~if~} z_h = 1
\end{array} 
\right.  \forall \, h \in \cH \label{nonc-oid} 
\end{align}
\end{subequations}
where the balance constraint~\eqref{nonc-balance} is enforced at each node $h \in
\cH$;~\eqref{nonc-oid}  indicates which inverters
have to be dispatched (i.e., $ (P_{c,h}, Q_{c,h})  \in
\cF^\mathrm{OID}_h$), or, operate in the business-as-usual mode (i.e., $ (P_{c,h},
Q_{c,h}) = (0,0)$); and, the constraint on $V_0$ is left implicit. Finally, 
$c_z \in \mathbb{R}_+$ is a weighting coefficient, used to trade off achievable cost $C(\bv, \bp_{c})$ 
for the number of controlled inverters. When $c_z$ represents a fixed reward for customers providing 
ancillary services~\cite{Caramanis_PESTD14} and $C(\bv, \bp_{c})$ models costs associated with active power losses and active power set points, OID~\eqref{Poidnonconvex} returns the inverter setpoints that minimize the economic cost
incurred by the feeder operation.

Unfortunately, problem~\eqref{Poidnonconvex} is nonconvex and it contains binary variables; 
thus, it is challenging to solve optimally and efficiently, 
even by utilizing off-the-shelf solvers for mixed-integer nonlinear programs. 
Nevertheless, a computationally-affordable \emph{convex} reformulation was introduced 
in~\cite{OID}, by leveraging contemporary sparsity-promoting regularization~\cite{Wiesel11} and
semidefinite relaxation (SDR) techniques~\cite{LavaeiLow, Dallanese-TSG13} as summarized next.

To address the non-convexity of constraints~\eqref{nonc-balance}--\eqref{nonc-Vlimits}, consider expressing 
powers and voltage magnitudes as linear functions of the
outer-product complex Hermitian matrix $\bV := \bv \bv^\sfH$,  and to
reformulate the OID problem with cost and constraints that are linear functions of
$\bV$. Specifically, define the matrix $\mathbf{Y}_n := \be_n \be_n^\sfT
\mathbf{Y}$ per node $n$, where $\{\mathbf{e}_{n}\}_{n \in \cN}$
denotes the canonical basis of $\mathbb{R}^{|\cN|}$. Further, based on
$\mathbf{Y}_n$, define also the Hermitian matrices $\bA_{n} :=
\frac{1}{2} (\bY_n + \bY_n^\sfH)$, $\bB_{n} := \frac{j}{2} (\bY_n -
\bY_n^\sfH) $, and $\bM_{n} := \be_n \be_n^\sfT$. Then, 
the node balance constraints for active and reactive powers can be equivalently 
expressed as $\trace(\bA_h \bV)   =  P_{h}^{\textrm{av}} - P_{c,h} - P_{\ell,h}$ and 
$\trace(\bB_h \bV)   = Q_{c,h} - Q_{\ell,h} $, respectively.  Similarly, constraint~\eqref{nonc-Vlimits}
can be equivalently expressed as $V_{\mathrm{min}}^2  \leq \trace(\bM_n \bV)  \leq V_{\mathrm{max}}^2$.
The technical constraints $\bV \succeq \mathbf{0}$ and
$\rank(\bV) = 1$ need to be added, to ensure recoverability of the voltage 
vector \bv~\cite{LavaeiLow, Dallanese-TSG13}. The
only source of non-nonconvexity is now constraint
$\rank(\bV) = 1$; however in the spirit of SDR, this constraint can be
dropped. If the optimal solution of the relaxed problem has
rank 1, then the resultant power flows are globally optimal for
given inverter setpoints. 

As for the binary variables $\{z_h\}$, notice first that if PV
inverter $h$ is \emph{not} selected for ancillary services, one 
has that $P_{c,h} = Q_{c,h} = 0$ [cf.~\eqref{nonc-oid}]. Thus,
assuming that only a subset of PV inverters may need to be controlled in order to 
ensure electrical network constraints and minimize~\eqref{nonc-cost}, 
one has that the $2 |\cH|\times 1$ real-valued vector
$[\bp_c^\sfT, \bq_c^\sfT ]^\sfT$ is
\emph{group sparse}~\cite{Wiesel11}; that is, either the $2
\times 1$ sub-vectors $[P_{c,h}, Q_{c,h}]^\sfT$ equal to $\mathbf{0}$,
or not. In lieu of binary
variables, this group-sparsity attribute enables PV inverter selection by regularizing the
cost in~\eqref{Poidnonconvex} with the following group-sparsity-promoting function:
\begin{align}
G(\bp_{c},\bq_{c})  :=  c_z \sum_{h \in \cH} \,  \|[P_{c,h}, Q_{c,h}]\|_2 .  \label{Glasso_powers} 
\end{align}

Leveraging these tools, a relaxation of the OID problem is obtained as:
\begin{subequations} 
\label{Poid}
\begin{align}
& \hspace{-1.9cm} \min_{\bV, \bp_{c}, \bq_{c}}   C(\bV, \bp_{c}) + G(\bp_{c},\bq_{c}) \label{Pm-cost}   \\
\,\, \mathrm{s.\,to} \,\,  \bV  \succeq \mathbf{0}, &  \mathrm{~and}  \nonumber  \\ 
\trace(\bA_h \bV)  & = - P_{\ell,h} + P_{h}^{\textrm{av}} - P_{c,h} \hspace{.65cm} \forall \, h \in \cH \label{Pm-balanceP}  \\
\trace(\bB_h \bV)  & = - Q_{\ell,h} + Q_{c,h}   \hspace{1.5cm} \forall \, h \in \cH \label{Pm-balanceQ} \\
\trace(\bA_n \bV)  &  = 0,~\trace(\bB_n \bV) = 0 \hspace{1.05cm} \forall \, n \in \cU  \label{Pm-balancePoles}  \\
 V_{\mathrm{min}}^2 & \leq \trace(\bM_n \bV)  \leq V_{\mathrm{max}}^2  \hspace{.85cm}  \forall \, n \in \cN   \label{Pm-Vlimits} \\
& \hspace{-.8cm} (P_{c,h}, Q_{c,h}) \in \cF^\mathrm{OID}_h(P_{h}^{\textrm{av}}) \hspace{1.05cm} \forall \, h \in \cH. \label{Pm-oid}  
\end{align}
\end{subequations}
Problem~\eqref{Poid} is convex, and can be readily re-formulated in a standard semidefinite programming (SDP) form by resorting to the epigraph forms of $G(\bp_{c},\bq_{c})$ and $C(\bV, \bp_{c})$~\cite{Vandenberghe96}, as well as the linear matrix inequality form of $Q_{c,h}^2  \leq  S_{h}^2 - (P_{h}^{\textrm{av}}  - P_{c,h})^2$ obtainable by using the Shur complement.  

When the distribution system is balanced and radial,  sufficient conditions for obtaining a rank-$1$ solution in SDR-based OPF-type reformulations are available in~\cite{Lavaei_tree14,Zhang12}, and they can be conveniently tailored to~\eqref{Poid}; for example, one requirement is that the cost~\eqref{Pm-cost} is increasing in the injected active powers.  
What is more, constraint $\bV  \succeq \mathbf{0}$ can be equivalently re-written as $\bV^{(i,j)}  \succeq \mathbf{0}, \forall \, (i,j) \in \cE$, with $\bV^{(i,j)}$ denoting the $2 \times 2$ sub-matrix of $\bV$ corresponding to nodes $i$ and $j$. Since $|V_n| > 0$ for all nodes, one has that each constraint $\bV^{(i,j)}  \succeq \mathbf{0}$ can be further re-expressed as  ($\bV_{ij}$ is the $(i,j)$-th entry of $\bV$)
\begin{align}
\bV_{ii} > 0 , \, \bV_{jj} > 0 , \, \bV_{ij} = \bV^*_{ji} , \textrm{~and~} |\bV_{ij}|^2 - \bV_{ii} \bV_{jj} \leq 0 ,  \nonumber 
\end{align}
which is a second-order cone constraint. Thus, for radial and balanced topologies,~\eqref{Poid} can be transformed into a second-order cone program, with due computational advantages. The worst-case complexity of an SDP is on the order $\mathcal{O}(N_v^{4.5} \log(1/\epsilon))$ for general-purpose solvers, with $N_v$ denoting the total number of variables in the problem, and $\epsilon > 0$ a given solution accuracy~\cite{Vandenberghe96,Nesterov94}. The worst-case complexity of a second-order cone program is on the order of $\mathcal{O}(N_v^3 \log(1/\epsilon))$~\cite{Nesterov94}. Notice however that sparsity in $\{\bA_{n},\bB_{n}, \bM_{n}\}$ and the chordal structure of the underlying electrical graph can be exploited to devise customized solvers with reduced computational burden; see e.g.,~\cite{Jabr12}. Finally, extensions of~\eqref{Poid} to multi-phase unbalanced distribution systems can be derived by following the method in~\cite{Dallanese-TSG13}. 

\vspace{0.1cm}

\noindent \emph{Remark (ZIP load model)}. A constant power load model is utilized in the OID~\eqref{Poidnonconvex}. However, the OID formulation can be broadened to account for constant-impedance, constant-current, and constant-power load components (i.e., the so-called ZIP model~\cite{Chassin11}) by following the method in~\cite{molzahn_zip14}. 
For constant impedance loads, the demands are proportional to the voltage magnitude squared; thus, they can be easily incorporated in~\eqref{Poid}. On the other hand, since constant current loads are functions of the voltage magnitudes, appropriate reformulations of~\eqref{Poid} are required. In particular,~\cite{molzahn_zip14} suggests to replace matrix $\bV$ with $\tilde{\bV}:= \tilde{\bv} \tilde{\bv}^\sfH$, where the voltage-related vector $\tilde{\bv}$ is defined as $\tilde{\bv} := [1, \bv^\sfT]^\sfT$ (clearly, matrices $\bA_n$, $\bB_n$, and $\bM_n$ are re-defined accordingly). In this case, it may not be possible to find a rank-$1$ matrix $\tilde{\bV}$, although the obtained solution yields a reasonably good approximation of the constant current loads; for more details, see the discussion provided in~\cite{molzahn_zip14}.

\section{Inverter Dispatch Under Forecasting Errors} \label{sec:Robust}

For given available powers $\{P_{h}^{\textrm{av}}\}_{h \in \cH}$, the OID task~\eqref{Poid} identifies the inverters that must provide ancillary services in order to avoid overvoltage conditions [cf.~\eqref{Pm-Vlimits}], and computes the steady-state setpoints that minimize the selected operational and economic objectives [cf.~\eqref{Pm-cost}]. This approach is suitable for real-time network operation, based on instantaneous measurements of loads and available powers. On the other hand, solar irradiance forecasting errors must be taken into account when OID is utilized for ancillary-services procurement, in either day-ahead or hour-ahead ancillary service markets~\cite{Craciun13,Caramanis_PESTD14}. In this case, system operators cannot rely on the expected value of the active power available from the PV array  to quantify the amount of  ancillary services to provision. In fact, an excess of generation (compared to the expected one) may require additional inverters other than the ones scheduled without accounting for forecasting errors, to deviate from the business-as-usual setpoint~\cite{Liu08,Cpuc12}.

In the remainder of this section, the so-called CVaR will be utilized to capture the risk of excess in the active power generation, and subsequently proactively select the inverters that will be required to provide ancillary services.

\subsection{Overview of the Conditional Value-at-Risk Approach}
\label{sec:CVaR}

An overview of the value-at-risk (VaR) and CVaR --- measures typically
considered in risk-aware portfolio optimization~\cite{Rockafellar00} --- 
is given in this subsection. These tools will be utilized to formulate 
the uncertainty-aware OID in Section~\ref{sec:robustOID}. 

Suppose $\bp^{\textrm{av}} := [P_{1}^{\textrm{av}}, \ldots,
P_{|\cH|}^{\textrm{av}}]^\sfT$ is a real-valued random vector, and let
$\rho(\bp^{\textrm{av}})$ denote its probability density
function. Assume that $\rho(\bp^{\textrm{av}})$ is known (or an
empirical estimate is available~\cite{Bacher2009,Lorenz09}), and
supported on a closed and bounded set $\cD \subset
\mathbb{R}^{|\cH|}$. For example, in the context of solar irradiance
forecasting, a viable choice for $\rho(\bp^{\textrm{av}})$ would be a
truncated multivariate Gaussian distribution, as described
in~\cite{Bacher2009}. See also e.g.,~\cite{Lorenz09} for additional
models for $\rho(\bp^{\textrm{av}})$ in the context of solar
irradiance forecasting. 

Let $r: \mathbb{R}^{|\cH|} \times \cD \rightarrow \mathbb{R}$ be a
real-valued function of both the random vector $\bp^{\textrm{av}}$ and
the \emph{vector of presumed powers} $\bd \in \mathbb{R}^{|\cH|}$. 
In particular,  let $r(\bd,\bp^{\textrm{av}}) =
\sum_{h\in\cH} [P_{h}^{\textrm{av}} - d_h]_+$ capture possible
excess of power generation during hours with high and yet uncertain
generation (and, hence, the risk of overvoltages throughout the
distribution feeder).\footnote{Another viable choice is $r(\bd,\bp^{\textrm{av}}) =
  [\sum_{h\in\cH} (P_{h}^{\textrm{av}} - d_h)]_+$; that is, the
  network-wide surplus of active power. However,
  $r(\bd,\bp^{\textrm{av}}) = \sum_{h\in\cH} [P_{h}^{\textrm{av}} -
  d_h]_+$ captures local (as opposed to network-wide) random changes
  in the active power injections, and it is therefore a more suitable
  indicator for the risk of high active power flows in sections of the
  feeder.} Henceforth, $r(\bd,\bp^{\textrm{av}})$ will be referred to as the 
\emph{surplus generation function}. Notice that
$r(\bd,\bp^{\textrm{av}})$ takes positive values only when
$P_{h}^{\textrm{av}} > d_h$ for at least one inverter. For a given
vector $\bd$, $r(\bd,\bp^{\textrm{av}})$ is a
random variable with cumulative distribution function
\begin{align}
\hspace{-.2cm} \Psi_r(\bd,\alpha) := \Prob\{r(\bd,\bp^{\textrm{av}}) \leq \alpha \} = \int_\cD \mathbb{I}_{\left\{r(\bd,\bp) \leq \alpha\right\}} \rho(\bp) \mathrm{d} \bp . \hspace{-.1cm}
\label{CVaRcdf} 
\end{align}
Notice that $\Psi_r(\bd,\alpha)$ is continuous from the right (but not
necessarily from the left), nondecreasing in $\alpha$, and
parameterized by $\bd$~\cite{Rockafellar00}. Intuitively, $
\Psi_r(\bd,\alpha)$ quantifies the probability of the actual available
power exceeding the presumed value $\bd$. Based on~\eqref{CVaRcdf}, the
VaR and CVaR measures are defined next (see~\cite{McNeilbook,Rockafellar00} for additional details).

\begin{figure}[t]
\begin{center}
\includegraphics[width=0.40\textwidth]{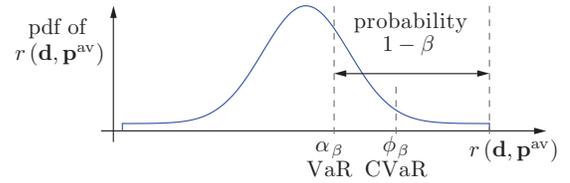}
\caption{Illustrative example of the CVaR associated with function $r(\bd,\bp^{\textrm{av}})$. }
\label{F_distribution}
\vspace{-.5cm}
\end{center}
\end{figure}

For a user-prescribed probability level $\beta \in (0,1)$, the
corresponding VaR, denoted as $\alpha_\beta$, associated with the
random surplus generation function $r(\bd,\bp^{\textrm{av}})$, is the
left endpoint of the non-empty interval collecting the values of
$\alpha$ for which $\Psi_r(\bd,\alpha) = \beta$; i.e.,
\begin{align}
\label{VaR} 
\alpha_\beta(\bd) :=  \inf \left\{ \alpha \in \mathbb{R} : \Psi_r(\bd,\alpha) \geq \beta \right\} .
\end{align}
For any $\bd$, the CVaR, denoted as $\phi_\beta(\bd)$, is the expected value of the surplus
generation function when considering entries that are
greater than or equal to $\alpha_\beta(\bd)$:
\begin{align}
\phi_\beta(\bd) := \frac{1}{1 - \beta} \int_\cD \mathbb{I}_{\left\{r(\bd,\bp) \geq \alpha_\beta(\bd) \right\}} r(\bd,\bp)  \rho(\bp) \mathrm{d} \bp.
\label{CVaR} 
\end{align}
In other words, in the $(1-\beta)$ percent of cases where
$r(\bd,\bp^{\textrm{av}}) = \sum_{h\in\cH} [P_{h}^{\textrm{av}} -
d_h]_+ \geq \alpha_\beta(\bd)$, the CVaR $\phi_\beta(\bd)$ quantifies the
expected amount of available active power further exceeding
$\alpha_\beta(\bd)$. An illustrative example of the VaR and CVaR associated with a random 
function $r(\bd,\bp^{\textrm{av}})$ is provided in Fig.~\ref{F_distribution}; 
in this example, the probability density function of $r(\bd,\bp^{\textrm{av}})$ is 
a truncated Gaussian. The CVaR is typically preferred
over VaR as a risk measure, since it is coherent (in fact, VaR
violates sub-additivity---one of the properties of a coherent
measure~\cite{McNeilbook}). 

Key to utilizing CVaR as a performance objective 
in risk-aware optimization tasks is the link established in~\cite{Rockafellar00} between $\phi_\beta(\bd)$ and the following real-valued function: 
\begin{align}
R_\beta(\alpha, \bd)  :=  \alpha + \frac{1}{(1-\beta)}  \int_\cD \left[ r(\bd,\bp)  - \alpha \right]_+ \rho(\bp) \mathrm{d} \bp  .\label{CVaRpowers} 
\end{align}
Specifically,~\cite[Thm.~1]{Rockafellar00} asserts the following three facts: 

\noindent $\mathrm{(F1)}$ $R_\beta(\alpha, \bd)$ is convex and continuously differentiable in $\alpha$.

\noindent $\mathrm{(F2)}$ For any $\bd   \in \mathbb{R}^{|\cH|}$, the CVaR $\phi_\beta(\bd)$ represents the minimum value of $R_\beta(\alpha, \bd)$; that is, 
\begin{align}
\phi_\beta(\bd) = \min_{\alpha \in \mathbb{R}} R_\beta(\alpha, \bd) .
\label{CVaRrel1} 
\end{align}

\noindent $\mathrm{(F3)}$ The set of minimizers 
\begin{align}
\cA_\beta(\bd) := \arg \min_{\alpha \in \mathbb{R}} R_\beta(\alpha, \bd)  
\label{CVaRrel2} 
\end{align}
is closed and bounded, and the VaR $\alpha_\beta(\bd)$ is the left endpoint of this interval.

An advantage of the integral function~\eqref{CVaRpowers} is that an
empirical estimate of $R_\beta(\alpha, \bd)$ can be obtained
via sample averaging. This is especially useful in cases when the
integral in~\eqref{CVaRpowers} cannot be evaluated in
closed-form. For instance, given $S$ Monte Carlo samples, $\{\bp^{\textrm{av}}[s] \in
\cD\}_{s = 1}^S$, of the random vector $\bp^{\textrm{av}}$, a
distribution-free approximation of $R_\beta(\alpha, \bd)$ is given by
\begin{align}
\hat{R}_\beta(\alpha, \bd)  =  \alpha + \frac{1}{S(1-\beta)} \sum_{s = 1}^{S} \left[  r(\bd,\bp^{\textrm{av}}[s])  - \alpha \right]_+,   \hspace{-.4cm} \label{CVaRapprox} 
\end{align}
and, for a sufficiently high number of samples $S$, almost sure
convergence of $\hat{R}_\beta(\alpha, \bd)$ to $R_\beta(\alpha, \bd)$
is guaranteed by the (strong) law of large numbers. Compared
to~\eqref{CVaRpowers}, the sample average $\hat{R}_\beta(\alpha, \bd)$
is not differentiable due to the projection operator
$[\cdot]_+$. However, this hurdle can be easily overcome by resorting
to the epigraph form of $\hat{R}_\beta(\alpha,
\bd)$~\cite{BoVa04}.

To consider $\hat{R}_\beta(\alpha, \bd)$ (or, its epigraph form) in multi-objective optimization problems, convexity with respect to (wrt) both $\alpha$ and $\bd$ is desirable. To this end, the following claims established in~\cite[Thm.~2]{Rockafellar00} can be conveniently leveraged: 

\noindent $\mathrm{(C1)}$ If function $r(\bd,\bp)$ is convex in $\bd$, then  $R_\beta(\alpha, \bd)$ is jointly convex in $\bd$ and $\alpha$; and  $\phi_\beta(\bd)$ is convex in $\bd$.     

\noindent $\mathrm{(C2)}$ The following equality holds
\begin{align}
\min_{\bd \in \mathbb{R}^{|\cH|}} \phi_\beta(\bd) =  \min_{\bd \in \mathbb{R}^{|\cH|}, \alpha \in \mathbb{R}} R_\beta(\alpha, \bd),  
\label{CVaRrel3} 
\end{align}
and $\bd^\star, \alpha^\star$ are minimizers of $R_\beta(\alpha, \bd)$ if and only if $\bd^\star$ is a minimizer of $\phi_\beta(\bd)$ and $ \alpha^\star \in \cA_\beta(\bd^\star)$.

Claim $\mathrm{(C2)}$ asserts that minimizing the CVaR wrt to the variables $\bd$ is equivalent to jointly minimizing $R_\beta(\alpha, \bd)$ (and thus $\hat{R}_\beta(\alpha, \bd)$) over $\bd$ and $\alpha$, with the VaR $\alpha^\star$ coming out as a byproduct. This feature will be exploited in the risk-aware OID framework outlined next, where $\bd$ represents the vector of presumed available powers associated with a given CVaR.

\subsection{Risk-aware inverter dispatch}
\label{sec:robustOID}

For $r(\bd,\bp^{\textrm{av}}) = \sum_{h\in\cH} [P_{h}^{\textrm{av}} - d_h]_+$,  function $R_\beta(\alpha, \bd)$ is jointly convex in $\bd$ and $\alpha$ by virtue of $(\mathrm{C1})$. Further, given $S$ independent samples $\{\bp^{\textrm{av}}[s] \in \cD\}_{s = 1}^S$, an approximation of $R_\beta(\alpha, \bd)$ is given by [cf.~\eqref{CVaRapprox}]
\begin{align}
\hat{R}_\beta(\alpha, \bd)  =  \alpha + \frac{1}{S(1-\beta)} \sum_{s = 1}^{S} \left[ \sum_{h\in\cH} \left[P_{h}^{\textrm{av}}[s] - d_h\right]_+  - \alpha \right]_+ . \hspace{-.4cm} \label{CVaRpowersApprox} 
\end{align}

Thus, given $\beta$ and the Monte Carlo samples $\{\bp^{\textrm{av}}[s] \in \cD\}_{s = 1}^S$, the objective of the risk-aware OID problem is to jointly minimize the OID objective~\eqref{Pm-cost} under the presumed available power levels $\bd$, as well as the risk of additional available power surplus, subject to the AC power flow and OID-related inverter constraints; that is,   
\begin{subequations} 
\label{Poid_robust}
\begin{align}
& \hspace{-1.9cm} \min_{\bV, \bp_{c}, \bq_{c}, \bd, \alpha}   C(\bV, \bp_{c},\bd) + G(\bp_{c},\bq_{c}) + c_R \hat{R}_\beta(\alpha, \bd)  \label{Pm-cost2}   \\
\,\, \mathrm{s.\,to} \,\,  \bV  \succeq \mathbf{0}, &  \mathrm{~and}  \nonumber  \\ 
\trace(\bA_h \bV)  & = - P_{\ell,h} + d_{h} - P_{c,h} \hspace{.8cm} \forall \, h \in \cH \label{Pm-balanceP2}  \\
\trace(\bB_h \bV)  & = - Q_{\ell,h} + Q_{c,h}   \hspace{1.5cm} \forall \, h \in \cH \label{Pm-balanceQ2} \\
\trace(\bA_n \bV)  &  = 0,~\trace(\bB_n \bV) = 0 \hspace{1.05cm} \forall \, n \in \cU  \label{Pm-balancePoles2}  \\
 V_{\mathrm{min}}^2 & \leq \trace(\bM_n \bV)  \leq V_{\mathrm{max}}^2  \hspace{.85cm}  \forall \, n \in \cN   \label{Pm-Vlimits2} \\
& \hspace{-.8cm} (P_{c,h}, Q_{c,h}) \in \cF^\mathrm{OID}_h(d_h) \hspace{1.2cm} \forall \, h \in \cH \label{Pm-oid2}  \\
 \bd & \in \cD \label{Pm-probSupport2}
\end{align}
\end{subequations}
where $c_R \in \mathbb{R}_+$ is a predetermined parameter, used to trade off achievable CVaR values for OID objectives at the $\beta$-risk level. Problem~\eqref{Poid_robust} is convex, and can be re-stated in either standard SDP form by using the epigraph form of~\eqref{Pm-cost2}~\cite{Vandenberghe96}, or, in standard SOCP form for systems that are radial and balanced.

To appreciate the usefulness of the CVaR risk measure, suppose that for given $\beta$,  it turns out that at least $H < |\cH|$ inverters are required to curtail at most $\{\bar{P}_{c,h}\}$ W, in order to ensure voltage regulation in the $\beta$ percent of the cases; that is, whenever $\sum_{h\in\cH} [P_{h}^{\textrm{av}} - d_h]_+ \leq \alpha_\beta(\bd)$. Then, minimizing the CVaR $\phi_\beta(\bd)$ is equivalent to minimizing the \emph{additional} amount of active powers that inverters may be required to curtail in case of unexpected over-generation (i.e., when $\sum_{h\in\cH} [P_{h}^{\textrm{av}} - d_h]_+ \geq \alpha_\beta(\bd)$), or, minimizing the number of \emph{additional} inverters that may be called upon to provide ancillary services.
Elaborating further on the impacts of uncertainties on the system operational costs, suppose that function $C(\bV, \bp_{c},\bd)$ is set to  $C(\bV, \bp_{c},\bd) = c_L ( \mathbf{1}_{|\cH|}^\sfT (\bd - \bp_{c}) + P_0) + c_P  \mathbf{1}_{|\cH|}^\sfT \bp_{c}$, where the first term captures the cost incurred by power losses in the network, the available powers are $\bd$, and the second term models the cost of active power that can be curtailed; see e.g.,~\cite{Caramanis_PESTD14}. Further, recall that $G(\bp_{c},\bq_{c})$ accounts for possible fixed rewards for customers when their inverters are called upon to provide ancillary reserves. If $c_R$ quantifies the economic loss incurred by overvoltages, then~\eqref{Pm-cost2} strikes a balance between system operational costs when operated at a risk level $\beta$, and the economic loss that the system may incur in case of unexpected generation surplus. Section~\ref{sec:Simulations} will elaborate further on how to trade off CVaR for the amount of ancillary services to be provisioned.

To re-state~\eqref{Poid_robust} in a standard SDP form (similar steps can be followed for the SOCP case), assume for simplicity that $C(\bV, \bp_{c},\bd)$ is linear in its arguments. Consider then introducing the non-negative auxiliary vector variable $\bz := [z_1, \ldots, z_{|\cH|}]^\sfT$, replace $G(\bp_{c},\bq_{c})$ with $c_z  \mathbf{1}_{|\cH|}^\sfT \bz$ in~\eqref{Pm-cost2}, and add constraints $ \|[P_{c,h}, Q_{c,h}]\|_2 \leq z_h$, for all $h \in \cH$.  Then, by introducing auxiliary variables  $\by \in \mathbb{R}^{S}$ and $\{\bu_s \in \mathbb{R}^{|\cH|}\}_{s = 1}^S$ to upper bound  the projection terms~\cite{Rockafellar00}, and by using the Schur complement to convert quadratic and conic constraints into linear inequality constraints~\cite{Vandenberghe96},~\eqref{Poid_robust} can be re-stated in the following standard SDP form:     
\begin{subequations} 
\label{Poid_robust2}
\begin{align}
& \hspace{-1.0cm} \min_{\substack{\bV, \bp_{c}, \bq_{c}, \bd \\ \alpha, \by \succeq \mathbf{0}, \bu_s  \succeq \mathbf{0} \\ \bz \succeq \mathbf{0}}}   C(\bV, \bp_{c},\bd) + c_z  \mathbf{1}_{|\cH|}^\sfT \bz + c_R \alpha + \frac{c_R}{S(1-\beta)} \mathbf{1}_{S}^\sfT \by     \label{Pm-costSDPref}   \\
\,\, \mathrm{s.\,to} \,\,   & \bV  \succeq \mathbf{0},~\eqref{Pm-balanceP2}-\eqref{Pm-Vlimits2},\mathrm{~and}  \nonumber  \\ 
& \hspace{-.2cm}  \left[ \begin{array}{ccc}
z_h & 0 & P_{c,h} \\
0 & z_h & Q_{c,h} \\
P_{c,h} & Q_{c,h} & z_h
\end{array}
\right] \succeq \mathbf{0} \hspace{1.7cm}  \forall \, h \in \cH
\label{SDPref-norm2} \\
& \hspace{-.2cm}  \left[ 
\begin{array}{ccc}
-S_h^2 & Q_{c,h} & d_h - P_{c,h} \\
Q_{c,h} & -1 & 0 \\
d_{h} - P_{c,h}  & 0 & -1
\end{array}
\right] \preceq \mathbf{0} \hspace{.2cm}  \forall \, h \in \cH  \hspace{-.2cm}
\label{SDPref-ConstReactive} \\
& \mathbf{0} \preceq \bp_{c}  \preceq  \bd \label{mg-SDPref-curt} \\
& \bq_{c} \preceq \tan \theta (\bd - \bp_{c}) \label{mg-SDPref-pf} \\
&  \hspace{-.4cm} - \bq_{c} \preceq \tan \theta (\bd - \bp_{c}) \label{mg-SDPref-pf2} \\
& \mathbf{1}_{|\cH|}^\sfT \bu_s  \leq \alpha + y_s, \hspace{1.6cm} \forall \, s = 1, \ldots, S  \label{Pm-SDPref-aux2} \\
&  \bp^{\mathrm{av}}[s] - \bd \preceq \bu_s, \hspace{1.7cm} \forall \, s = 1, \ldots, S \label{Pm-SDPref-aux3} \\
& \bd \in \cD  . \label{Pm-probSupport3} 
\vspace{-.3cm}
\end{align}
\end{subequations}
  
\noindent \emph{Remark (load uncertainty)}. Although this section focused on solar irradiance forecasting errors, uncertainty in active and reactive household demands can also be accounted for in the risk-aware OID framework. For example, for the active power demand, function $r(\bd,\bp^{\textrm{av}}) = \sum_{h\in\cH} [P_{h}^{\textrm{av}} - P_{\ell,h} - (d_h - \ell_h)]_+$ can be utilized to capture surplus of \emph{net} generated active power throughout the feeder, where both $P_{h}^{\textrm{av}}$ and $P_{\ell,h}$ are now random variables, and $\ell_h$ is the counterpart of $d_h$ for the demanded active power. Then, the risk-aware OID problem is obtained by replacing $P_{\ell,h}$ with $\ell_h$ in~\eqref{Pm-balanceP2}.  A similar procedure can be followed for uncertain reactive loads.

\noindent \emph{Remark (optimal solution)}. On par with~\cite{Lavaei_tree14,Zhang12}, for distribution feeders that are radial and balanced, the semidefinite relaxation~\eqref{Poid_robust} is exact when the following sufficient conditions are satisfied: s1) the cost function is increasing with respect to the net active power injection; s2) the voltage angle difference $\theta_{ik}$ between nodes $i$ and $k$ is such that $-\tan^{-1}(b_{ik}/g_{ik}) \leq \theta_{ik} \leq \tan^{-1}(b_{ik}/g_{ik})$, with $y_{ik} = g_{ik} + j b_{ik}$ the admittance of the line $(i,j) \in \cE$; and, s3) inverters are able to absorb a ``sufficient'' amount of reactive power, with specific bounds quantified in~\cite[Thm.~1]{Zhang12}. Condition s2) is typically satisfied in practice, since voltage angle differences are small; condition s3) can be checked by inspecting $\cF^\mathrm{OID}$; and, s1) can be satisfied by appropriate tuning of the problem parameters.  For unbalanced feeders as well as meshed networks, efforts for finding sufficient conditions that ensure exactness of the semidefinite relaxation are still undergoing. However, the virtues of semidefinite relaxation have been demonstrated in e.g.,~\cite{Dallanese-TSG13} and~\cite{Gan14}.

\noindent \emph{Remark (uncertainty in the
    temperature)} The formulation~\eqref{Poid_robust} accounts for
  uncertainty in the available active power through function
  $R_\beta(\alpha, \bd)$.  Accordingly, it would be straightforward to
  translate forecasted irradiance \emph{and} temperature values into
  forecasted power values using standard models for PV modules and inverters, see, e.g.,~\cite{masters2004renewable}, provided the probability density
  function of irradiance and temperature forecasting errors are
  available [c.f.~\eqref{CVaRpowers} and~\eqref{CVaRpowersApprox}].
 
\section{Case Studies} \label{sec:Simulations}

To solve the OID problem, the distribution network operator requires: i) the Monte Carlo samples $\{\bp^{\textrm{av}}[s] \in \cD\}_{s = 1}^S$, based on the distribution of the solar irradiation error~\cite{Bacher2009}; ii) the network admittance matrix $\bY$; ii) the probability level $\beta$; the ratings $\{S_h\}$; and, v) the weighting coefficients $c_L, c_P, c_z, c_R$, which may be driven by ancillary service market strategies~\cite{Caramanis_PESTD14,Craciun13} and/or security-oriented objectives. The optimization package \texttt{CVX}\footnote{[Online] \texttt{http://cvxr.com/cvx/}} is  employed to solve the OID problem in \texttt{MATLAB}. In all the presented tests, the rank of matrix $\bV$ was always $1$, implying that the SDR relaxation for the power flow equations is tight.

\begin{figure}[t]
\begin{center}
\includegraphics[width=0.40\textwidth]{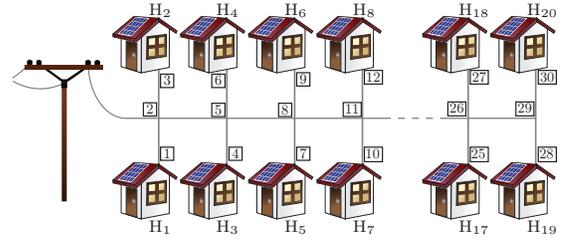}
\caption{Low-voltage residential network adopted for the case studies.}
\label{F_LV_network}
\vspace{-.5cm}
\end{center}
\end{figure}

The distribution network in Fig.~\ref{F_LV_network} is considered in the test cases, which is a larger version of the fishbone system utilized in~\cite{Tonkoski11,OID} to assess the impact of high PV generation in residential setups. The pole-pole distance is set to $30$ m, while the lengths of the drop lines are set to $20$ m. The values of the line impedances are adopted from~\cite{Tonkoski11}.

The 20 houses shown in Fig.~\ref{F_LV_network} feature fixed roof-top PV systems, with a DC-AC derating coefficient of $0.77$. The DC ratings of the houses are as follows: $5.52$ kW for houses $\mathrm{H}_1, \mathrm{H}_3, \mathrm{H}_{6}, \mathrm{H}_{7}, \mathrm{H}_{8}, \mathrm{H}_{9}, \mathrm{H}_{11}, \mathrm{H}_{14}, \mathrm{H}_{16}$, and $\mathrm{H}_{19}$; $8.00$ kW for houses $\mathrm{H}_2, \mathrm{H}_{10}, \mathrm{H}_{12}, \mathrm{H}_{13}, \mathrm{H}_{18}, \mathrm{H}_{20}$; and, $5.70$ kW for the remaining houses. The minimum power factor for the inverters is set to 0.85, and it is assumed that the PV inverters are oversized by $10\%$ of their AC rating~\cite{Turitsyn11}. To account for forecasting errors, the available powers are modeled as $P_h^{\textrm{av}} = \bar{P}_h^{\textrm{av}} + \Delta_{h}$, with $\bar{P}_h^{\textrm{av}}$ the (known) \emph{forecasted value} and $\Delta_{h}$ the (random) \emph{forecasting error}. The hourly forecasted values of the available powers $\{\bar{P}_h^{\textrm{av}}\}$ are computed using the System Advisor Model\footnote{[Online] \texttt{https://sam.nrel.gov/}.} of the National Renewable Energy Laboratory, based on  typical meteorological year data for Minneapolis, MN, during the month of July. Hourly PV generation in the interval $\cT := \{\mathrm{6AM}, \ldots, \mathrm{8PM}\}$  is considered. A zero-mean truncated Gaussian distribution is adopted for $\Delta_{h}$, with truncation at the $0.3$th and $99.7$th percentiles; see e.g.,~\cite{Bacher2009}. Random variables $\{\Delta_{h}\}$ are correlated across houses, and an exponentially decreasing correlation function $\E\{\Delta_{h} \Delta_{h^\prime} \} = \sigma_h \sigma_{h^\prime} e^{- d(h, h^\prime)/\tau}$  is used, where $\sigma_h$ is the standard deviation of $\Delta_{h}$, $d(h, h^\prime)$ the distance between houses $h$ and $h^\prime$, and $\tau = 300$ [m].  

The residential load profile is obtained from the Open Energy Info database,\footnote{[Online] \texttt{http://en.openei.org/datasets/node/961}} and the ``base load'' experienced in downtown Saint Paul, MN, during the month of July is used for this test case. To generate different load profiles, the base active power profile is perturbed using a truncated Gaussian random variable with zero mean and standard deviation $200$ W, and a power factor of 0.9 is presumed~\cite{Tonkoski11}. Finally, voltages $V^\textrm{min}$ and $V^\textrm{max}$ are set to 0.917 pu and 1.042 pu, respectively in this case study (see e.g., page 11 of the CAN3-C235-83 standard). 
The voltage at the secondary of the transformer is set at 1.02 pu, to ensure a minimum voltage magnitude of $0.917$ pu when the PV inverters do not generate power. 

\begin{figure}[t]
\centering
  \subfigure[]{\includegraphics[width=0.47\textwidth]{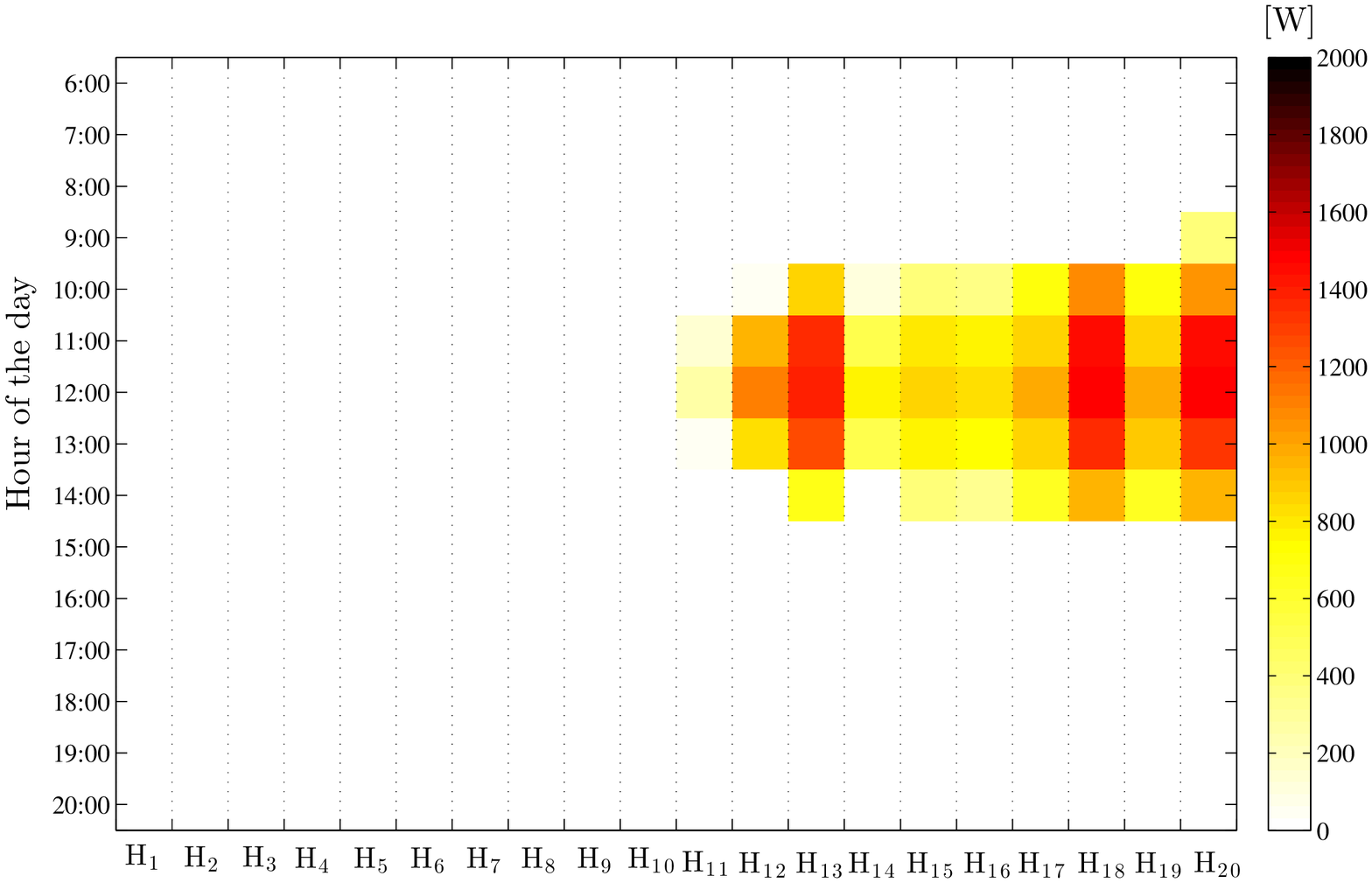}} 
  \subfigure[]{\includegraphics[width=0.47\textwidth]{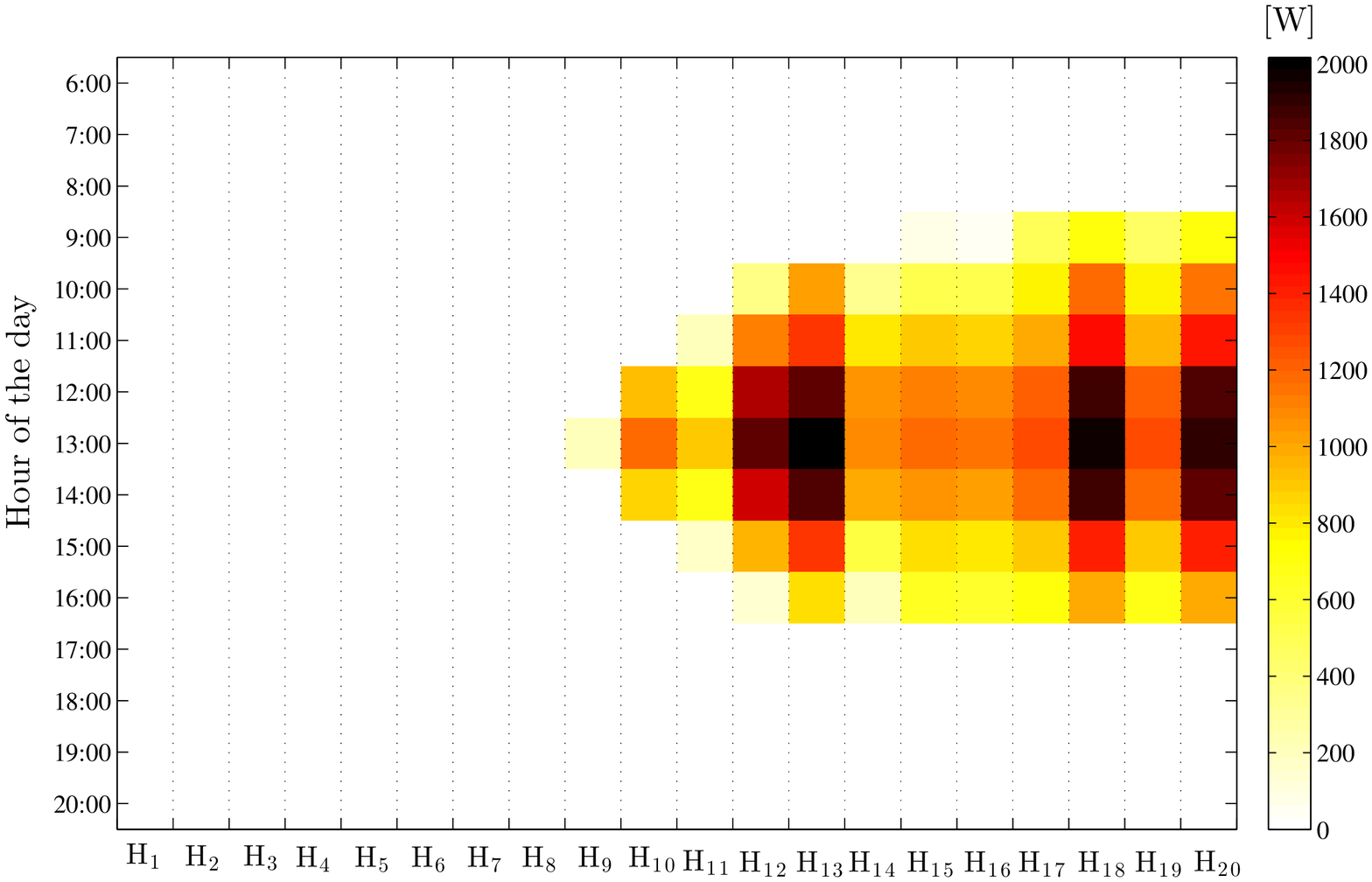}} 
  \caption{Dispatched inverters: provisioned active power curtailment for (a) $c_R = 0.01$ and (b) $c_R = 10$.} 
\label{F_OID_result} 
\end{figure}

\begin{table}[t]
\caption{Provisioned ancillary services for different risk levels ($\beta = 0.95, c_L = 1, c_P = 0.5, c_z = 0.9$)}
\vspace{-.2cm}
\begin{center}
\begin{tabular}{l||c|c|c|c}
 &  $P_c^{\mathrm{tot}}$ [kWh] & $Q_c^{\mathrm{tot}}$  [kVAr] & $N^{\mathrm{tot}}$ \\
\hline \hline
No risk & 27.83 & 4.65 &  44\\
$c_R = 0.01$ & 38.60  & 8.45 &  44 \\
$c_R = 0.1$   & 40.26  & 8.53 &  56 \\
$c_R = 1$      & 44.05  & 9.75 & 59 \\
$c_R = 10$   & 79.15 & 14.83 & 86 \\
\end{tabular}
\end{center}
\label{tab:test1_results}
\end{table}%

\begin{table}[t]
\caption{Provisioned ancillary services for different probability levels $\beta$ ($c_L = 1, c_P = 1, c_z = 0.9, c_R = 1$)}
\vspace{-.2cm}
\begin{center}
\begin{tabular}{l||c|c|c|c}
 &  $P_c^{\mathrm{tot}}$ [kWh] & $Q_c^{\mathrm{tot}}$  [kVAr] & $N^{\mathrm{tot}}$ \\
\hline \hline
No risk & 27.83 & 4.65 &  44\\
$\beta = 0.85$  & 41.05 & 9.51 &  55 \\
$\beta = 0.90$      & 41.35  & 9.65 &  57 \\
$\beta = 0.95$   &  44.05 & 9.75 & 59 \\
$\beta = 0.99$ & 49.55  & 10.22 &  64 \\
\end{tabular}
\end{center}
\label{tab:test2_results}
\end{table}%

\begin{table}[t]
\caption{Provisioned ancillary services for different uncertainty levels ($c_L = 1, c_P = 1, c_z = 0.9, c_R = 1, \beta = 0.95$)}
\vspace{-.2cm}
\begin{center}
\begin{tabular}{l||c|c|c|c}
$\sigma_h/\bar{P}_{h}^{\mathrm{av}}$ &  $P_c^{\mathrm{tot}}$ [kWh] & $Q_c^{\mathrm{tot}}$  [kVAr] & $N^{\mathrm{tot}}$ \\
\hline \hline
0 & 27.83 & 4.65 &  44\\
0.05 & 42.12 & 9.19 & 57 \\
0.10 &  44.05 & 9.75 & 59 \\
0.15 & 46.25 & 10.35 & 63 \\
0.20 & 51.93 & 11.64 & 65 \\
\end{tabular}
\end{center}
\label{tab:test3_results}
\end{table}%

In the first setup, the standard deviation of the solar power prediction error $\sigma_h$ amounts to $10 \%$ of the forecasted value~\cite{Bacher2009} (that is, $\sigma_h/\bar{P}_{h}^{\mathrm{av}} = 0.1$); $S = 1000$; $\beta = 0.95$; $c_z = 0.9$ to capture fixed rewards when inverters are called upon providing ancillary services;  cost $C(\bV, \bp_{c},\bd)$ is set to $C(\bV, \bp_{c},\bd) = c_L ( \mathbf{1}_{|\cH|}^\sfT (\bd - \bp_{c} - \bp_{\ell} ) + P_0) + c_P \mathbf{1}_{|\cH|}^\sfT \bp_{c}$, with $c_L = 1$ and $c_P = 0.5$. The amount of active power that can be curtailed by each inverter during the day is illustrated in Fig.~\ref{F_OID_result}, for $c_R = 0.01$ (lower weight given to the CVaR) and $c_R = 10$ (low CVaR objectives). 
It can be clearly seen that the amount of active power provisioned from each inverter increases with the increasing of $c_R$, thus ensuring an enhanced system protection against unexpected boosts in the solar irradiation. Clearly, the enhanced system protection comes at the expense of a higher reward for customers providing this ancillary service (modeled by the term $c_P \mathbf{1}_{|\cH|}^\sfT \bp_{c}$). As observed also in~\cite{OID}, inverters with higher ratings may be required to curtail more active power. To facilitate fairness among customers, the term  $\| \bPi \bp_c \|_2$ can be included in~\eqref{Pm-cost2}, where $\bPi := \bI_{|\cH|} - \frac{1}{|\cH|}\mathbf{1}_{|\cH| \times 1} \mathbf{1}_{|\cH| \times 1} ^\sfT$.
Finally, notice that for $c_R = 10$, an increased number of inverters are required to curtail active power. 

This trend is confirmed by the results reported in Table~\ref{tab:test1_results}, where: $P_c^{\mathrm{tot}} := \sum_{t \in \cT} \sum_{h \in \cH} P_{c,h}(t)$ denotes the conglomerate active power curtailment that is procured during the day, with $P_{c,h}(t)$ the amount of power that can be curtailed from inverter $h$ at time $t$;  $Q_c^{\mathrm{tot}} := \sum_{t \in \cT} \sum_{h \in \cH} |Q_{c,h}(t)|$ the overall reactive power procured for reactive support purposes;  and  $N^{\mathrm{tot}}$  the total number of inverters called upon providing ancillary services over $\cT$ (out of $20 \times |\cT| = 360$). These values are compared with the ``no risk'' setup, where the solar forecasting errors are neglected, and  ancillary services are provisioned by solving~\eqref{Poid} with  $\{P_{h}^{\mathrm{av}}\}$ replaced by $\{\bar{P}_{h}^{\mathrm{av}}\}$. Clearly,  considering only the expected available powers $\{\bar{P}_{h}^{\mathrm{av}}\}$ yields an underestimate of the amount of active and reactive power reserves that  may be needed to ensure voltage regulation. This corroborates the ability of the proposed approach to trade off risks of overvoltage conditions for the amount of ancillary services to be secured~\cite{Caramanis_PESTD14}.  

Table~\ref{tab:test2_results} quantifies the procured active and reactive reserves for different  values for the probability level $\beta$. In this setup, the other problem parameters are set as $c_L = 1, c_P = 1, c_z = 0.9$ and $c_R = 1$. With the increasing of $\beta$, progressively higher solar irradiation conditions are considered in the risk-aware OID problem [cf.~\eqref{VaR}]; this explains why the amount of active and reactive powers reserves that are secured by the OID framework increases.  A similar trend can be noticed in Table~\ref{tab:test3_results}, where values for the standard deviation of the forecasting errors are tested. 
Notice that the results for $\sigma_h = 0$ (i.e., perfect knowledge of the solar irradiation conditions) coincide with the ``no risk'' setup of Tables~\ref{tab:test1_results} and~\ref{tab:test2_results}. As expected, the higher is $\sigma_h$, the higher is the number of inverters that may be required to provide ancillary services.  

\section{Concluding Remarks} \label{sec:Conclusions}
 
The present paper dealt with ancillary service provisioning in distribution systems 
in the presence of solar irradiance forecasting errors.
The proposed uncertainty-aware OID identifies the subset of critical inverters that 
strongly impact both voltage profile and network performance objectives, and quantifies 
the amount of active and reactive powers to be procured from each inverter. 
The CVaR measure was utilized to capture (and minimize) the risk of overvoltages 
throughout the feeder. Although the formulated OID task involves the solution of a 
nonconvex mixed-integer nonlinear program, a convex relaxation was developed by leveraging  
sparsity-promoting regularization approaches and semidefinite relaxation techniques. 
Using real-world PV-generation and load-profile data, it is shown how the proposed framework can 
trade of the risk of PV generation surplus for the amount of ancillary services to be provisioned.

\appendix

The dependence between active power injections and voltage magnitudes in low-voltage 
distribution systems is briefly analyzed in the following.  

Consider a low-voltage single-phase distribution line, and let $Z := R + j \omega L$ be its impedance, where  $R$, $L$, and $\omega$ denote the per-unit-length line resistance and inductance, and electrical frequency. Typical values for $R$ and $L$ are on the order of  $10^{-1}$~$\Omega/\textrm{km}$  and $10^{-4}$~$\textrm{H}/\textrm{km}$, respectively (see e.g., the specifications of cables NS 90 3/0 AWG utilized for the pole-to-pole connections and cables NS 90 1/0 AWG for drop lines~\cite{Tonkoski11}). Let $Y := 1/Z = $, $G:= \Re\{Y\} = R/|Z|^2$ and $B:= \Im\{Y\} = - \omega L /|Z|^2$. Further, let $|V_1| e^{j \theta_1}$ and $|V_2| e^{j \theta_2}$ denote the phasors for the voltages at the two end points of the line. Although a single line is considered for simplicity, claims  naturally extend to low-voltage systems with arbitrary topologies. 

The active- and reactive-power injections at node 1 of the line are given by
\small
\begin{subequations}
\begin{align}
&P_1 =  |V_1|^2 G - |V_1| |V_2| G \cos(\theta) - |V_1| |V_2| B \sin (\theta), \label{RealPower1} \\
&Q_1 = -|V_1|^2 B + |V_1| |V_2| B \cos(\theta) - |V_1| |V_2| G \cos(\theta) ,  \label{ReactivePower1} 
\end{align}
\end{subequations}
\normalsize
where $\theta := \theta_1 - \theta_2$. Similarly, the active- and
reactive-power injections at node 2 of the line are given by 
\small
\begin{subequations}
\begin{align}
& P_2  =  |V_2|^2 G - |V_1| |V_2| G \cos(\theta) + |V_1| |V_2| B \sin (\theta) , \label{RealPower2} \\
& Q_2  = -|V_2|^2 B + |V_1| |V_2| B \cos(\theta) + |V_1| |V_2| G \cos(\theta). \label{ReactivePower2} 
\end{align}
\end{subequations}
\normalsize

\begin{figure*}[t!]
\normalsize
\begin{align}
\left[
\begin{array}{c}
\Delta P_1 \\
\Delta P_2 \\
\Delta Q_1 \\
\Delta Q_2 
\end{array}
\right] 
= 
\left[
\begin{array}{cccc}
2 G |V_1| - G |V_2| & - G |V_1|  & - B |V_1| |V_2|  & B |V_1| |V_2|   \\
- G |V_2| & - G |V_1| + 2 G |V_2|  & B | V_1| |V_2| & - B |V_1| |V_2| \\
-2  B |V_1| +  B |V_2| & B |V_1| & - G |V_1| |V_2| & G |V_1| |V_2| \\
B |V_2| & B |V_1|  -2 B |V_2| & G |V_1| |V_2|  & -B |V_1| |V_2|  
\end{array}
\right] 
\left[
\begin{array}{c}
\Delta |V_1| \\
\Delta |V_2| \\
\Delta \theta_1 \\
\Delta \theta_2 
\end{array}
\right]  . 
\label{sensitivity_matrix} 
\end{align}
\hrulefill
\end{figure*}

With regard to~\eqref{RealPower1}--\eqref{ReactivePower2}, define the sensitivity matrix 
\begin{align}
\bS(|V_1|, |V_2|, \theta_1, \theta_2) := 
\left[ 
\begin{array}{cccc}
\frac{\partial P_1}{\partial |V_1|} & \frac{\partial P_1}{\partial |V_2|} & \frac{\partial P_1}{\partial \theta_1} & \frac{\partial P_1}{\partial \theta_2} \\ 
\frac{\partial P_2}{\partial |V_1|} & \frac{\partial P_2}{\partial |V_2|} & \frac{\partial P_2}{\partial \theta_1} & \frac{\partial P_2}{\partial \theta_2} \\  
\frac{\partial Q_1}{\partial |V_1|} & \frac{\partial Q_1}{\partial |V_2|} & \frac{\partial Q_1}{\partial \theta_1} & \frac{\partial Q_1}{\partial \theta_2} \\  
\frac{\partial Q_2}{\partial |V_1|} & \frac{\partial Q_2}{\partial |V_2|} & \frac{\partial Q_2}{\partial \theta_1} & \frac{\partial Q_2}{\partial \theta_2} 
\end{array}
\right] , 
\nonumber
\end{align}
which relates power variations with perturbations of the voltage phasors around a given operational point. 

Next, assume small voltage angle variations; that is, $\theta \ll 1$, $\cos(\theta) \approx 1$ and $\sin(\theta) \approx \theta$. Under these assumptions, we can relate sensitivities of voltage magnitudes and angles to real and reactive power injections through~\eqref{sensitivity_matrix}. Furthermore, since $\frac{\omega L}{R} \ll 1$ in low-voltage distribution systems (this condition is not necessarily true in medium-voltage networks), it follows that $B \ll 1$, and thus the effects of voltage magnitude and phase variations on the complex powers approximately  decouples as
 \begin{align}
\left[ \hspace{-.2cm}
\begin{array}{c}
\Delta P_1 \\
\Delta P_2 
\end{array}
\hspace{-.2cm}
\right] 
\approx
\left[
\begin{array}{cccc}
2 G |V_1| - G |V_2| & - G |V_1|     \\
- G |V_2| & - G |V_1| + 2 G |V_2|   \\
\end{array}
\right] 
\left[ \hspace{-.2cm}
\begin{array}{c}
\Delta |V_1| \\
\Delta |V_2| 
\end{array}
\hspace{-.2cm}
\right], \nonumber
\end{align}
\vspace{-.4cm}
\begin{align}
\left[ \hspace{-.2cm}
\begin{array}{c}
\Delta Q_1 \\
\Delta Q_2 
\end{array}
\hspace{-.2cm}
\right] 
\approx
\left[
\begin{array}{cccc}
- G |V_1| |V_2| & G |V_1| |V_2| \\
 G |V_1| |V_2|  & -G |V_1| |V_2|  
\end{array}
\right] 
\left[ \hspace{-.2cm}
\begin{array}{c}
\Delta \theta_1 \\
\Delta \theta_2 
\end{array}
\hspace{-.2cm}
\right] .
\nonumber
\end{align}
Thus, due to the high resistance-to-inductance ratio in low-voltage distribution lines, voltage magnitudes are more sensitive to variations in the active power flows. It follows that curtailing active power during peak generation hours represents a viable way to avoid overvoltage conditions throughout the feeder. Furthermore, the higher is the solar irradiation, the higher is the overall amount of active power that should be curtailed in order to enforce voltage regulation.

\bibliographystyle{IEEEtran}
\bibliography{biblio}

\begin{thebibliography}{10}
\providecommand{\url}[1]{#1}
\csname url@samestyle\endcsname
\providecommand{\newblock}{\relax}
\providecommand{\bibinfo}[2]{#2}
\providecommand{\BIBentrySTDinterwordspacing}{\spaceskip=0pt\relax}
\providecommand{\BIBentryALTinterwordstretchfactor}{4}
\providecommand{\BIBentryALTinterwordspacing}{\spaceskip=\fontdimen2\font plus
\BIBentryALTinterwordstretchfactor\fontdimen3\font minus
  \fontdimen4\font\relax}
\providecommand{\BIBforeignlanguage}[2]{{%
\expandafter\ifx\csname l@#1\endcsname\relax
\typeout{** WARNING: IEEEtran.bst: No hyphenation pattern has been}%
\typeout{** loaded for the language `#1'. Using the pattern for}%
\typeout{** the default language instead.}%
\else
\language=\csname l@#1\endcsname
\fi
#2}}
\providecommand{\BIBdecl}{\relax}
\BIBdecl

\bibitem{Liu08}
E.~Liu and J.~Bebic, ``Distribution system voltage performance analysis for
  high-penetration photovoltaics,'' Feb. 2008, {NREL} Technical Monitor: B.
  Kroposki. Subcontract Report {NREL/SR}-581-42298.

\bibitem{Cpuc12}
{California Public Utilities Commission}, ``Advanced inverter technologies
  report,'' Jan. 2013, [Online] \texttt{http://www.cpuc.ca.gov}.

\bibitem{Tonkoski11}
R.~Tonkoski, L.~A.~C. Lopes, and T.~H.~M. El-Fouly, ``Coordinated active power
  curtailment of grid connected {PV} inverters for overvoltage prevention,''
  \emph{IEEE Trans. on Sust. Energy}, vol.~2, no.~2, pp. 139--147, Apr. 2011.

\bibitem{Tricoli10}
A.~F. Vizoso, L.~Piegari, and P.~Tricoli, ``A photovoltaic power unit providing
  ancillary services for smart distribution networks,'' in \emph{Intl. Conf. on
  Renewable Energies and Power Quality}, Las Palmas, Spain, 2010.

\bibitem{Caramanis_PESTD14}
E.~Ntakou and M.~C. Caramanis, ``Price discovery in dynamic power markets with
  low-voltage distribution-network participants,'' in \emph{IEEE PES Trans. \&
  Distr. Conf.}, Chicago, IL, 2014.

\bibitem{Craciun13}
D.~Craciun and D.~Geibel, ``Evaluation of ancillary services provision
  capabilities from distributed energy supply,'' in \emph{Intl. Conf. on
  Electricity Distribution {CIRED}}, Stockholm, June 2012.

\bibitem{Turitsyn11}
K.~Turitsyn, P.~Sulc, S.~Backhaus, and M.~Chertkov, ``Options for control of
  reactive power by distributed photovoltaic generators,'' \emph{Proc. of the
  IEEE}, vol.~99, no.~6, pp. 1063--1073, 2011.

\bibitem{Aliprantis13}
P.~Jahangiri and D.~C. Aliprantis, ``Distributed {Volt/VAr} control by {PV}
  inverters,'' \emph{IEEE Trans. Power Syst.}, vol.~28, no.~3, pp. 3429--3439,
  Aug. 2013.

\bibitem{Bolognani13}
S.~Bolognani and S.~Zampieri, ``A distributed control strategy for reactive
  power compensation in smart microgrids,'' \emph{IEEE Trans. on Autom.
  Control}, vol.~58, no.~11, pp. 2818--2833, 2013.

\bibitem{Samadi14}
A.~Samadi, R.~Eriksson, L.~Soder, B.~G. Rawn, and J.~C. Boemer, ``Coordinated
  active power-dependent voltage regulation in distribution grids with pv
  systems,'' \emph{IEEE Trans. on Power Del.}, vol.~29, no.~3, pp. 1454--1464,
  June 2014.

\bibitem{OID}
E.~Dall'Anese, S.~V. Dhople, and G.~B. Giannakis, ``Optimal dispatch of
  photovoltaic inverters in residential distribution systems,'' \emph{IEEE
  Trans. Sustainable Energy}, vol.~5, no.~2, pp. 487--497, Apr. 2014.

\bibitem{OID_TEC}
E.~Dall'Anese, S.~V. Dhople, B.~B. Johnson, and G.~B. Giannakis,
  ``Decentralized optimal dispatch of photovoltaic inverters in residential
  distribution systems,'' \emph{IEEE Trans. on Energy Conversion}, 2014, to
  appear. Also available at: \texttt{http://arxiv.org/abs/1403.1341}.

\bibitem{McNeilbook}
A.~McNeil, R.~Frey, and P.~Embrechts, \emph{Quantitative Risk Management:
  Concepts Techniques and Tools}.\hskip 1em plus 0.5em minus 0.4em\relax
  Princeton University Press, 2005.

\bibitem{Rockafellar00}
R.~T. Rockafellar and S.~Uryasev, ``Optimization of conditional
  value-at-risk,'' \emph{Journal of Risk}, vol.~2, no.~3, pp. 21--41, 2000.

\bibitem{Bacher2009}
P.~Bacher, H.~Madsen, and H.~A. Nielsen, ``Online short-term solar power
  forecasting,'' \emph{Solar Energy}, vol.~83, no.~10, pp. 1772--1783, 2009.

\bibitem{Lorenz09}
E.~Lorenz, J.~Hurka, D.~Heinemann, and H.~G. Beyer, ``Irradiance forecasting
  for the power prediction of grid-connected photovoltaic systems,'' \emph{IEEE
  J. Sel. Top. App. Earth Observ. and Rem. Sensing}, vol.~2, no.~1, pp. 2--10,
  2009.

\bibitem{Sjodin12}
E.~Sjodin, D.~F. Gayme, and U.~Topcu, ``Risk-mitigated optimal power flow for
  wind powered grids,'' in \emph{American Control Conf.}, Montreal, Canada,
  June 2012.

\bibitem{Warrington13}
J.~Warrington, P.~Goulart, S.~Mariethoz, and M.~Morari, ``Policy-based reserves
  for power systems,'' \emph{IEEE Trans. Power Syst.}, vol.~28, no.~4, pp.
  4427--4437, Nov 2013.

\bibitem{Summers14}
T.~Summers, J.~Warrington, M.~Morari, and J.~Lygeros, ``Stochastic optimal
  power flow based on convex approximations of chance constraints,'' in
  \emph{Power Systems Computation Conf.}, Wroclaw, Poland, Aug. 2014.

\bibitem{Wang12}
Q.~Wang, Y.~Guan, and J.~Wang, ``A chance-constrained two-stage stochastic
  program for unit commitment with uncertain wind power output,'' \emph{IEEE
  Trans. Power Syst.}, vol.~27, no.~1, pp. 206--215, Feb 2012.

\bibitem{Chertkov13}
D.~Bienstock, M.~Chertkov, and S.~Harnett, ``Chance constrained optimal power
  flow: Risk-aware network control under uncertainty,'' 2013, [Online]
  \texttt{http://arxiv.org/abs/1209.5779}.

\bibitem{Kloppel13}
M.~Kloppel, A.~Gabash, A.~Geletu, and P.~Li, ``Chance constrained optimal power
  flow with non-gaussian distributed uncertain wind power generation,'' in
  \emph{Intl. Conf. on Environment and Electrical Engineering}, May 2013, pp.
  265--270.

\bibitem{YuSGComm13}
Y.~Zhang and G.~B. Giannakis, ``Robust optimal power flow with wind integration
  using conditional value-at-risk,'' in \emph{4th Intl. Conf. on Smart Grid
  Communications}, Vancouver, Canada, Oct. 2013.

\bibitem{Dallanese-TSG13}
E.~Dall'Anese, H.~Zhu, and G.~B. Giannakis, ``Distributed optimal power flow
  for smart microgrids,'' \emph{IEEE Trans. Smart Grid}, vol.~4, no.~3, pp.
  1464--1475, Sep. 2013.

\bibitem{LavaeiLow}
J.~Lavaei and S.~H. Low, ``Zero duality gap in optimal power flow problem,''
  \emph{IEEE Trans. Power Syst.}, vol.~1, no.~1, pp. 92--107, 2012.

\bibitem{kerstingbook}
W.~H. Kersting, \emph{Distribution System Modeling and Analysis}.\hskip 1em
  plus 0.5em minus 0.4em\relax 2nd ed., Boca Raton, {FL}: {CRC} Press, 2007.

\bibitem{Wiesel11}
A.~T. Puig, A.~Wiesel, G.~Fleury, and A.~O. Hero, ``Multidimensional
  shrinkage-thresholding operator and group {LASSO} penalties,'' \emph{IEEE
  Sig. Proc. Letters}, vol.~18, no.~6, pp. 363--366, Jun. 2011.

\bibitem{Vandenberghe96}
L.~Vandenberghe and S.~Boyd, ``Semidefinite programming,'' \emph{{SIAM}
  Review}, vol.~38, no.~1, pp. 49--95, Mar. 1996.

\bibitem{Lavaei_tree14}
J.~Lavaei, D.~Tse, and B.~Zhang, ``Geometry of power flows and optimization in
  distribution networks,'' \emph{IEEE Trans. Power Syst.}, vol.~29, no.~2, pp.
  572--583, March 2014.

\bibitem{Zhang12}
B.~Zhang, A.~Y. Lam, A.~Dominguez-Garcia, and D.~Tse, ``An optimal and
  distributed method for voltage regulation in power distribution systems,''
  \emph{IEEE Trans. Power Syst.}, to appear. See also:
  \texttt{http://arxiv.org/abs/1204.5226}.

\bibitem{Nesterov94}
A.~Nesterov and Y.~Nemirovski, \emph{Interior-Point Polynomial Algorithms in
  Convex Programming}.\hskip 1em plus 0.5em minus 0.4em\relax SIAM Studies in
  Applied Mathematics, 1994.

\bibitem{Jabr12}
R.~A. Jabr, ``Exploiting sparsity in {SDP} relaxations of the {OPF} problem,''
  \emph{IEEE Trans. Power Syst.}, vol.~2, no.~27, pp. 1138--1139, May 2012.

\bibitem{Chassin11}
F.~Chassin, E.~Mayhorn, M.~Elizondo, and S.~Lus, ``Load modeling and
  calibration techniques for power system studiess,'' in \emph{North American
  Power Symp.}, Boston, MA, Aug. 2011.

\bibitem{molzahn_zip14}
D.~K. Molzahn, B.~C. Lesieutre, and C.~L. DeMarco, ``Approximate representation
  of {ZIP} loads in a semidefinite relaxation of the {OPF} problem,''
  \emph{IEEE Trans. Power Syst., Letters}, vol.~29, no.~4, pp. 1864--1865, July
  2014.

\bibitem{BoVa04}
S.~Boyd and L.~Vandenberghe, \emph{Convex Optimization}.\hskip 1em plus 0.5em
  minus 0.4em\relax Cambridge University Press, 2004.

\bibitem{Gan14}
L.~Gan and S.~H. Low, ``Convex relaxations and linear approximation for optimal
  power flow in multiphase radial networks,'' 2014, [Online]
  \texttt{http://arxiv.org/abs/1406.3054}.

\bibitem{masters2004renewable}
G.~Masters, \emph{Renewable and Efficient Electric Power Systems}.\hskip 1em
  plus 0.5em minus 0.4em\relax Wiley, 2004.

\end{thebibliography}

\end{document}